\documentclass[10pt, letterpaper]{amsart}
\usepackage[latin1]{inputenc}
\usepackage{amsfonts}
\usepackage{amsmath}
\usepackage{amssymb}
\usepackage{amsthm}
\usepackage[all]{xy}
\usepackage{graphicx}
\usepackage[margin=1.25in]{geometry}

\theoremstyle{plain}
	\newtheorem*{theorem*}{Theorem}
	\newtheorem{theorem}{Theorem}[section]
	\newtheorem*{lemma*}{Lemma}
	\newtheorem{lemma}{Lemma}[section]
	\newtheorem*{claim*}{Claim}
	
	\newtheorem*{prop*}{Proposition}
	\newtheorem{prop}{Proposition}[section]
	
	\newtheorem*{cor*}{Corollary}

\theoremstyle{definition}
	\newtheorem{note}{Remark}[section]
	\newtheorem{definition}{Definition}[section]
	
\newcommand{\Z}{\mathbb{Z}}
\newcommand{\R}{\mathbb{R}}

\newcommand{\C}{\mathbb{C}}

\newcommand{\noqed}{\renewcommand{\qedsymbol}{}}

\begin{document}

\title{Fundamental solution for $(\Delta - \lambda_z)^{\nu}$ on a symmetric space $G/K$}
\author{Amy T. DeCelles}
\address{University of St. Thomas, Department of Mathematics, 2115 Summit Avenue, St. Paul, MN 55105-1079}
\email{adecelles@stthomas.edu, amy.decelles@gmail.com}
\urladdr{http://cam.mathlab.stthomas.edu/decelles}
\thanks{This paper presents results from the author's PhD thesis, completed under the supervision of Professor Paul Garrett, whom the author thanks warmly.  The author would also like to thank Brian Hall for several helpful conversations.  The author was partially supported by the Doctoral Dissertation Fellowship from the Graduate School of the University of Minnesota and by NSF grant DMS-0652488.}

\begin{abstract}We determine a fundamental solution for the differential operator $(\Delta - \lambda_z)^{\nu}$ on the Riemannian symmetric space $G/K$,  where $G$ is any complex semi-simple Lie group, and $K$ is a maximal compact subgroup.  We develop a global zonal spherical Sobolev theory, which enables us to use the harmonic analysis of spherical functions to obtain an integral representation for the solution.  Then we obtain an explicit expression for the fundamental solution, which allows relatively easy estimation of its behavior in the eigenvalue parameter $\lambda_z$, with an eye towards further applications to automorphic forms involving associated Poincar\'{e} series.\end{abstract}

\maketitle

\section{Introduction}
\label{free_sp_fund_soln_intro}

We determine a fundamental solution for the differential operator $(\Delta - \lambda_z)^{\nu}$ on the Riemannian symmetric space $G/K$, where $G$ is any complex semi-simple Lie group and $K$ is a maximal compact subgroup.  Since the delta distribution $\delta_{1\cdot K}$ at the base point is a \emph{bi-$K$-invariant}, \emph{compactly supported} distribution, a suitable global zonal spherical Sobolev theory ensures that the harmonic analysis of spherical functions produces a solution.  In this paper, we first develop the suitable Sobolev theory, then derive the fundamental solution.    To our knowledge, this is the first construction of Sobolev spaces of bi-$K$-invariant compactly supported distributions. 


Instead of using the existence of a fundamental solution to prove solvability of a differential operator, as in, for example,  \cite{andersen1998, andersen2000,ban-schlichtkrull1993,benabdallah-rouviere1984,bopp-harinck1992}, we obtain an \emph{explicit expression} for the fundamental solution, with eye towards further applications involving the associated Poincar\'{e} series.   For example, we have already obtained an explicit formula relating the number of lattice points in an expanding region in a symmetric space to the automorphic spectrum \cite{decelles-lattice2011}.

In particular, the presence of a complex (eigenvalue) parameter $z$ in the differential  operator makes the fundamental solution suitable for further applications, and the simple, explicit nature of the fundamental solution allows relatively easy estimation of its behavior in the eigenvalue parameter, proving convergence of the associated Poincar\'{e} series in $L^2$ and, in fact, in a Sobolev space sufficient to prove continuity \cite{decelles-lattice2011}.  Further, this makes it possible to determine the vertical growth of the Poincar\'{e} series in the eigenvalue parameter.

For a derivation of the fundamental solution in the case $G = SL_2(\C)$, \emph{assuming} a suitable global zonal spherical Sobolev theory, see \cite{garrett-newark2010,garrett2010-harm}.  Our results for the general case are sketched in \cite{garrett-durham2010}.  After having submitted an initial version of this paper, it was brought to our attention that Wallach derives a similar, though less explicit, formula in Section 4 of \cite{wallach}.  An introduction to positively indexed Sobolev spaces of bi-$K$-invariant functions can be found in \cite{anker1992}.



Our main result is the following theorem, whose proof is given in \ref{free_sp_fund_soln_pvthm}.

\begin{theorem*}  Let $G$ be a complex semi-simple Lie group with maximal compact $K$.  When $G$ is of odd rank, let $\nu = d + \frac{n+1}{2}$, where $d$ is the number of positive roots, not counting multiplicities, and $n= \mathrm{dim} \,\mathfrak{a}$ the rank.  Then the bi-$K$-invariant fundamental solution $u_z$ for the operator $(\Delta - \lambda_z)^{\nu}$ on $G/K$ is given by:
$$u_z(a) \;\; = \;\;  \frac{(-1)^{d+(n+1)/2}  \, \pi^{(n+1)/2}}{\pi^+(\rho)\, \Gamma(d+(n+1)/2) } \; \cdot \; \prod_{\alpha\in \Sigma^+}  \, \frac{\alpha(\log a)}{2 \sinh(\alpha(\log a))} \; \cdot \;\; \frac{e^{-z|\log a|}}{z}  $$
When $G$ is of even  rank, let $\nu  =  d + \frac{n}{2} + 1$.  Then, with $K_n$ the usual modified Bessel function,
$$u_z(a) \;\; = \;\; \frac{(-1)^{d + (n/2) + 1} \,\pi^{n/2}}{\pi^+(\rho) \; \Gamma(d + (n/2) + 1)} \; \cdot  \; \prod_{\alpha\in \Sigma^+}  \, \frac{\alpha(\log a)}{2 \sinh(\alpha(\log a))} \;\cdot \; \frac{ |\log a|}{ z} \; \cdot \; K_1(z \, |\log a|)$$
\end{theorem*}


\section{Spherical transforms, global zonal spherical Sobolev spaces, \\and differential equations on $G/K$}

\subsection{Spherical transform and inversion}

Let $G$ be a complex semi-simple Lie group with finite center and $K$ a maximal compact subgroup.  Let $G = NAK$, $\mathfrak{g} = \mathfrak{n} + \mathfrak{a} + \mathfrak{k}$ be corresponding Iwasawa decompositions.  Let $\Sigma$ denote the set of roots of $\mathfrak{g}$ with respect to $\mathfrak{a}$, let $\Sigma^+$ denote the subset of positive roots (for the ordering corresponding to $\mathfrak{n}$), and let $\rho = \tfrac{1}{2} \sum_{\alpha \in \Sigma^+} m_{\alpha} \alpha$, $m_{\alpha}$ denoting the multiplicity of $\alpha$.  Let $\mathfrak{a}_{\C}^{\ast}$ denote the set of complex-valued linear functions on $\mathfrak{a}$.  Let $X = K \backslash G/K$ and $\Xi = \mathfrak{a}^{\ast}/W \approx \mathfrak{a}_+$.  
The spherical transform of Harish-Chandra and Berezin integrates a bi-$K$-invariant against a zonal spherical function:
$$\mathcal{F}f \, (\xi) \; = \; \int_G f(g) \, \overline{\varphi}_{\rho + i\xi}(g) \, dg $$
Zonal spherical functions $\varphi_{\rho + i \xi}$ are eigenfunctions for Casimir (restricted to bi-$K$-invariant functions) with eigenvalue $\lambda_{\xi} = -(|\xi|^2 + |\rho|^2)$.   The inverse transform is
$$ \mathcal{F}^{-1} f \; = \; \int_{\Xi} f(\xi) \, \varphi_{\rho + i\xi} \; |\mathbf{c}(\xi)|^{-2} d\xi$$
where $\mathbf{c}(\xi)$ is the Harish-Chandra $\mathbf{c}$-function and $d\xi$ is the usual Lebesgue measure on $\mathfrak{a}^{\ast} \approx \R^n$.  For brevity, denote $L^2(\Xi, |\mathbf{c}(\xi)|^{-2})$ by $L^2(\Xi)$.  The Plancherel theorem asserts that the spectral transform and its inverse are isometries between $L^2(X)$ and $L^2(\Xi)$.  

\subsection{Characterizations of Sobolev spaces}

We define positive index zonal spherical Sobolev spaces as left $K$-invariant subspaces of completions of $C_c^{\infty}(G/K)$ with respect to a topology induced by seminorms associated to derivatives from the universal enveloping algebra, as follows.  Let $\mathcal{U}\mathfrak{g}^{\leq \ell}$ be the finite dimensional subspace of the universal enveloping algebra $\mathcal{U}\mathfrak{g}$ consisting of elements of degree less than or equal to $\ell$.  Each $\alpha \in \mathcal{U}\mathfrak{g}$ gives a seminorm  $\nu_{\alpha}(f) \; = \; \lVert \alpha f \rVert_{L^2(G/K)}^2$  on $C_c^{\infty}(G/K)$.

\begin{definition}Consider the space of smooth functions that are bounded with respect to these seminorms:
$$\{ f \in C^{\infty}(G/K) : \nu_{\alpha} f \, < \, \infty\;\text{ for all } \alpha \in \mathcal{U}\mathfrak{g}^{\leq \ell}\}$$
Let $H^{\ell}(G/K)$ be the completion of this space with respect to the topology induced by the family $\{\nu_{\alpha} : \alpha \in \mathcal{U}\mathfrak{g}^{\leq \ell} \}$.  The \emph{global zonal spherical Sobolev space} $H^{\ell}(X) = H^{\ell}(G/K)^K$ is the subspace of left-$K$-invariant functions in $H^{\ell}(G/K)$.\end{definition}

\begin{prop} The space of test functions $C_c^{\infty}(X)$ is dense in $H^{\ell}(X)$.


\begin{proof}  We approximate a smooth function $f \in H^{\ell}(X)$ by pointwise products with smooth cut-off functions, whose construction (given by \cite{gangolli-varadarajan}, Lemma 6.1.7) is as follows.  Let $\sigma(g)$ be the geodesic distance between the cosets $1 \cdot K$ and $g \cdot K$ in $G/K$.    For $R>0$, let $B_R$ denote the ball $B_R \; = \; \{ g \in G:  \sigma(g) < R \}$.  Let $\eta$ be a non-negative smooth bi-$K$-invariant function, supported in $B_{1/4}$, such that $\eta(g) = \eta(g^{-1})$, for all $g \in G$.  Let $\mathrm{char}_{R+1/2}$ denote the characteristic function of $B_{R + 1/2}$, and let $\eta_R \; = \; \eta \, \ast \, \mathrm{char}_{R+1/2} \; \ast \; \eta$.  As shown in \cite{gangolli-varadarajan}, $\eta_R$ is smooth, bi-$K$-invariant, takes values between zero and one, is identically one on $B_R$ and identically zero outside $B_{R+1}$, and, for any $\gamma \in \mathcal{U}\mathfrak{g}$, there is a constant $C_\gamma$ such that
$$\sup_{g \in G} \; | (\gamma \, \eta_R) (g) | \;\; \leq \;\; C_\gamma$$

We will show that the pointwise products $\eta_R \cdot f$ approach $f$ in the $\ell^{\mathrm{th}}$ Sobolev topology, i.e. for any $\gamma \in \mathcal{U}\mathfrak{g}^{\leq \ell}$, $ \nu_{\gamma} \big(\eta_R \cdot f - f\big) \to 0$ as $R \to \infty$.  By definition,
$$ \nu_{\gamma} \big(\eta_R \cdot f - f\big) \;\; = \;\; \lVert \gamma \big( \eta_R \cdot f - f\big) \rVert_{L^2(G/K)}$$
Leibnitz' rule implies that $\gamma\big( \eta_R \cdot f - f\big)$ is a finite linear combination of terms of the form $\alpha(\eta_R - 1) \cdot \beta f$ where $\alpha$, $\beta \in \mathcal{U}\mathfrak{g}^{\leq \ell}$.  When $\mathrm{deg}(\alpha) = 0$, 
$$\lVert \alpha(\eta_R -1)\cdot \beta f \rVert_{L^2(G/K)}\;\; \ll \;\; \lVert (\eta_R -1) \cdot \beta f \rVert_{L^2(G/K)} \;\; \leq \;\;  \int_{\sigma(g) \, \geq \, R} \; |(\beta \, f) (g)|^2 \, dg$$
Otherwise, $\alpha (\eta_R -1) = \alpha \eta_R$, and
\begin{eqnarray*}
\lefteqn{\lVert \alpha(\eta_R -1)\cdot \beta f \rVert_{L^2(G/K)} \;\; = \;\; \lVert \alpha \eta_R \cdot \beta f \rVert_{L^2(G/K)}}\\
&  \ll &  \sup_{g \in G} \; |\alpha \, \eta_R(g) | \; \cdot  \; \int_{\sigma(g) \, \geq \, R} \; |(\beta \, f) (g)|^2 \, dg \;\;  \ll \;\; \int_{ \sigma(g)\, \geq \, R} \; |(\beta \, f) (g)|^2 \, dg
\end{eqnarray*}
Let $B$ be any bounded set containing all of the (finitely many) $\beta$ that appear as a result of applying Leibniz' rule.  Then
$$ \nu_{\gamma} \big(\eta_R \cdot f - f\big) \;\; \ll \;\; \sup_{\beta \in B} \;\;  \int_{ \sigma(g) \, \geq \, R} \; |(\beta \, f) (g)|^2 \, dg$$
Since $B$ is bounded and $f \in H^{\ell}(X)$, the right hand side approaches zero as $R \to \infty$.
\end{proof}
\end{prop}

\begin{prop} \label{casimir_norm} Let $\Omega$ be the Casimir operator in the center of $\mathcal{U}\mathfrak{g}$.  The norm $\lVert \,\cdot \, \rVert_{2\ell}$ on $C_c^{\infty}(G/K)^K$ given by
$$\lVert f \rVert_{2\ell}^2 \;\; = \;\; \lVert f \rVert^2 \; + \; \lVert (1-\Omega) \, f \rVert^2 \; + \; \lVert (1-\Omega)^2 \, f \rVert^2 \; + \; \dots \; + \; \lVert (1-\Omega)^{\ell} \, f \rVert^2$$
where $\lVert \, \cdot \, \rVert$ is the usual norm on $L^2(G/K)$, induces a topology on $C^{\infty}_c(G/K)^{K}$ that is equivalent to the topology induced by the family $\{\nu_{\alpha}: \alpha \in \mathcal{U}\mathfrak{g}^{\leq \, 2\ell}\}$ of seminorms and with respect to which $H^{2\ell}(X)$ is a Hilbert space.

\begin{proof}Let $\{X_i\}$ be a basis for $\mathfrak{g}$ subordinate to the Cartan decomposition $\mathfrak{g} = \mathfrak{p} + \mathfrak{k}$.  Then $\Omega \; = \; \sum_{i } X_i \, X_i^{\ast}$, where $\{X_i^\ast\}$ denotes the dual basis, with respect to the Killing form.  Let $\Omega_{\mathfrak{p}}$ and $\Omega_{\mathfrak{k}}$ denote the subsums corresponding to $\mathfrak{p}$ and $\mathfrak{k}$ respectively.  Then $\Omega_{\mathfrak{p}}$ is a non-positive operator, while $\Omega_\mathfrak{k}$ is non-negative.
\begin{lemma} For any non-negative integer $r$, let $\Sigma_{r}$ denote the finite set of possible $K$-types of $\gamma \, f$, for $\gamma \in \mathcal{U}\mathfrak{g}^{\leq r}$ and $f \in C_c^{\infty}(G/K)^K$, and let $C_r$ be a constant greater than all of the finitely many eigenvalues $\lambda_\sigma$ for $\Omega_\mathfrak{k}$ on the $K$-types $\sigma \in \Sigma_r$.  For any $\varphi \in C_c^{\infty}(G/K)$ of $K$-type $\sigma \in \Sigma_m$ and $\beta = x_1 \dots x_n$ a monomial in $\mathcal{U}\mathfrak{g}$ with $x_i \in \mathfrak{p}$,
$$\langle \beta \, \varphi, \beta \, \varphi \rangle \;\;  \leq \;\; \langle (-\Omega + C_{m+n-1})^n \, \varphi, \, \varphi \rangle$$
where $\langle \, , \, \rangle$ is the usual inner product on $L^2(G/K)$.
\begin{proof}
We proceed by induction on $n = \mathrm{deg} \, \beta$.  For $n=1$, $\beta = x \in \mathfrak{p}$.
Let  $\{X_i\}$  be a self-dual basis for $\mathfrak{p}$ such that $X_1 = x$.  Then, 
$$\langle x   \varphi, \, x  \varphi \rangle \;\; \leq \;\; \sum_i \langle X_i \,  \varphi, \, X_i \, \varphi \rangle \;\; = \;\;   - \sum_i \langle X_i^2 \,  \varphi, \varphi \rangle  \;\; = \;\;  \langle -\Omega_{\mathfrak{p}} \, \varphi,  \varphi \rangle \;\; = \;\;   \langle (-\Omega + \Omega_{\mathfrak{k}}) \,  \varphi,  \varphi \rangle $$
$$\;\;\;\;\;  \;\; \leq \;\;   \langle (-\Omega + C_m) \, \varphi,  \varphi \rangle \;\; = \;\;   \langle (-\Omega + C_{m+n-1}) \, \varphi,  \varphi \rangle$$
For $n > 1$, write $\beta = x \gamma$, where $x = x_1$ and  $\gamma = x_2 \dots x_n$.  Then the $K$-type of $\gamma \varphi$ lies in $\Sigma_{m + n-1}$, and by the above argument,
$$\langle x \, \gamma \varphi, \; x \, \gamma \varphi \rangle \;\; \leq \;\; \langle (-\Omega + C_{m+n-1})\, \gamma \varphi, \, \gamma \varphi \rangle$$
Let  $C_c^{\infty}(G/K)_{\Sigma_{r}}$ be the subspace of $C_c^{\infty}(G/K)$ consisting of functions of $K$-type in $\Sigma_{r}$ and $L^2(G/K)_{\Sigma_{r}}$ be the corresponding subspace of $L^2(G/K)$.  For the moment, let $\Sigma = \Sigma_{m+n-1}$ and  $C = C_{m+n-1}$. Then, by construction,  $-\Omega_\mathfrak{k} + C$ is positive on $C_c^{\infty}(G/K)_{\Sigma}$, and thus $-\Omega + C \, = \, -\Omega_\mathfrak{p} - \Omega_{\mathfrak{k}} + C$ is a positive densely defined symmetric operator on $L^2(G/K)_{\Sigma}$.  Thus, by Friedrichs \cite{friedrichs34, friedrichs35}, there is an everywhere defined inverse $R$, which is a positive symmetric \emph{bounded} operator on $L^2(G/K)_{\Sigma}$, and which, by the spectral theory for bounded symmetric operators, has a positive symmetric square root $\sqrt{R}$ in the closure of the polynomial algebra $\C[R]$ in the Banach space of bounded operators on $L^2(G/K)_\Sigma$.  Thus $-\Omega + C$ has a symmetric positive square root, namely $\big(\sqrt{R}\big)^{-1}$, defined on $C_c^\infty(G/K)_\Sigma$, commuting with all elements of $\mathcal{U}\mathfrak{g}$, and 
$$ \langle (-\Omega + C)\, \gamma \varphi, \, \gamma \varphi \rangle \;\; = \;\;  \langle \gamma \; \sqrt{-\Omega + C} \;  \varphi,  \;  \gamma \; \sqrt{-\Omega + C} \;  \varphi\rangle$$
Now the $K$-type of $\sqrt{-\Omega + C} \, \varphi$, being the same as that of $\varphi$, lies in $\Sigma_m$, so by inductive hypothesis,
\begin{eqnarray*} \langle \gamma \; \sqrt{-\Omega + C} \;  \varphi,  \;  \gamma \; \sqrt{-\Omega + C} \;  \varphi\rangle & \leq & \langle (-\Omega + C_{m+n-2})^{n-1} \; \sqrt{-\Omega + C} \;  \varphi,  \;  \sqrt{-\Omega + C} \;  \varphi\rangle\\
& = &  \langle (-\Omega + C_{m+n-2})^{n-1}(-\Omega + C_{m+n-1}) \,  \varphi,  \;   \varphi\rangle\\
& \leq &  \langle (-\Omega + C_{m+n-1})^{n} \,  \varphi,  \;   \varphi\rangle
\end{eqnarray*}
and this completes the proof of the lemma.
\noqed \end{proof}
\end{lemma}

Let $\alpha \in \mathcal{U}\mathfrak{g}^{\leq 2\ell}$.   By the Poincar\'{e}-Birkhoff-Witt theorem we may assume $\alpha$ is a monomial of the form $\alpha \; = \; x_1 \dots x_n \; y_1 \dots y_m$ where $x_i \in \mathfrak{p}$ and $y_i \in \mathfrak{k}$.  Then, for any $f \in C_c^{\infty}(G/K)^K$,
$$\nu_{\alpha} f \;\; = \;\; \langle \alpha f , \alpha f \rangle_{L^2(G/K)} \;\; = \;\; \langle x_1 \dots x_n \; f, \; x_1 \dots x_n f \rangle_{L^2(G/K)} \;\;\;\;\; (x_i \in \mathfrak{p})$$
By the lemma, there is a constant $C$, depending on the degree of $\alpha$, such that  $\nu_{\alpha}(f) \, \ll \, \langle (-\Omega + C)^{\mathrm{deg} \, \alpha} \, f , f \rangle$ for all $ f \in C_c^{\infty}(G/K)^K$.  In fact, for bi-$K$-invariant functions, $(-\Omega + C)^{\mathrm{deg} \, \alpha} \, f \; = \; (-\Omega_{\mathfrak{p}} + C)^{\mathrm{deg} \, \alpha} \, f $.  Since $\Omega_{\mathfrak{p}}$ is positive semi-definite, multiplying by a positive constant does not change the topology.  Thus, we may take $C=1$.  That is, the subfamily $\{\nu_{\alpha} : \alpha = (1-\Omega)^{k}, k \leq \ell \}$ of seminorms on $C_c^{\infty}(G/K)^{K}$ dominates the family $\{\nu_{\alpha} : \alpha \in \mathcal{U}\mathfrak{g}^{\leq \, 2\ell} \}$ and thus induces an equivalent topology.  
\end{proof}
\end{prop}

It will be necessary to have another description of Sobolev spaces.  Let
$$W^{2, \ell}(G/K) \;\; = \;\; \{f \in L^2(G/K): \alpha \, f \in L^2(G/K) \, \text{ for all } \alpha \in \mathcal{U}\mathfrak{g}^{\leq \ell} \}$$
where the action of $\mathcal{U}\mathfrak{g}$ on $L^2(G/K)$ is by distributional differentiation.  Give $W^{2,\ell}(G/K)$ the topology induced by the seminorms $\nu_{\alpha} f \; = \; \lVert \alpha \, f \rVert_{L^2(G/K)}^2$, $ \alpha \in \mathcal{U}\mathfrak{g}^{\leq \ell}$.  Let $W^{2, \ell}(X)$ be the subspace of left $K$-invariants.

\begin{prop} \label{large_small} These spaces are equal to the corresponding Sobolev spaces: 
$$W^{2, \ell}(G/K) \; = \; H^{\ell}(G/K) \;\;\;\;\; \text{ and } \;\;\;\;\; W^{2, \ell}(X)\; = \;H^{\ell}(X)$$

\begin{proof}
It suffices to show the density of test functions in $W^{2, \ell}(G/K)$.  Since $G$ acts continuously on $W^{2,\ell}(G/K)$ by left translation,  \emph{mollifications} are dense in $W^{2,\ell}(G/K)$; see \ref{gpints_moll}.  By Urysohn's Lemma, it suffices to consider mollifications of continuous, compactly supported functions.  Let $\eta \in C_c^{\infty}(G)$ and $f \in C_c^0(G/K)$.  Then, $\eta \cdot f$ is a smooth vector, and for all $\alpha \in \mathcal{U}\mathfrak{g}$, $\alpha \cdot (\eta \cdot f) \; = \; (L_{\alpha} \eta) \cdot f$.  For $X \in \mathfrak{g}$, the action on $\eta \cdot f$ as a vector is
$$X \cdot (\eta \cdot f) \;\; = \;\; \frac{\partial}{\partial t}\bigg|_{t=0} e^{tX} \; \cdot \;\; \int_G \eta(g) \;\; g\cdot f \, dg \;\; = \;\; \frac{\partial}{\partial t}\bigg|_{t=0}  \;\; \int_G \eta(g) \;\; (e^{tX} g)\cdot f \, dg $$
Now using the fact that $f$ is a \emph{function} and the group action on $f$ is by \emph{translation},
$$\big(X\cdot (\eta \cdot f)\big)(h) \;\; = \;\;  \frac{\partial}{\partial t}\bigg|_{t=0}  \;\; \int_G \eta(g) \; f(g^{-1}e^{-tX}h) \, dg \;\; = \;\; \frac{\partial}{\partial t}\bigg|_{t=0}  \;\; (\eta \cdot f)(e^{-tX}h)$$
Thus the smoothness of $(\eta \cdot f)$ as a vector implies that it is a genuine smooth function.  The support of $\eta \cdot f$ is contained in the product of the compact supports of $\eta$ and $f$.  Since the product of two compact sets is again compact, $\eta \cdot f $ is compactly supported.
\end{proof}
\end{prop}

\begin{note}
By Proposition \ref{casimir_norm}, $H^{2\ell}(X) = W^{2,2\ell}(X)$ is a Hilbert space with norm
$$\lVert f \rVert^2_{2\ell} \;\; = \;\; \lVert f \rVert^2 \; + \; \lVert (1-\Omega) \, f \rVert^2 \; + \; \dots \; + \; \lVert (1-\Omega)^{\ell} \, f \rVert^2$$
where $\lVert \, \cdot \, \rVert$ is the usual norm on $L^2(G/K)$, and $(1-\Omega)^k \, f$ is a distributional derivative.
\end{note}

\subsection{Spherical transforms and differentiation on Sobolev spaces}

\begin{prop}\label{diff_components} For $\ell \geq 0$, the Laplacian extends to a continuous linear map $H^{2\ell +2 }(X) \to H^{2\ell}(X)$; the spherical transform extends to a map on $H^{2\ell}(X)$; and
$$\mathcal{F} \big( (1-\Delta) f \big) \;\; = \;\; (1-\lambda_{\xi}) \cdot \mathcal{F}f \;\;\;\;\; \text{for all } f \in H^{2\ell+2}(X)$$
\begin{proof}
By the construction of the Sobolev topology, the Laplacian is a continuous map 
$$\Delta: C^{\infty}(G/K) \cap H^{2\ell+2}(G/K) \to C^{\infty}(G/K) \cap H^{2\ell}(G/K)$$
Since the Laplacian preserves bi-$K$-invariance, it extends to a (continuous linear) map, also denoted $\Delta$, from $H^{2\ell +2}(X)$ to $H^{2\ell}(X)$.  The spherical transform, defined on $C^{\infty}_c(G/K)^K$ by the integral transform of Harish-Chandra and Berezin, extends by continuity to $H^{2\ell}(X)$.  This extension agrees with the extension to $L^2(X)$ coming from Plancherel.  By integration by parts, $\mathcal{F}\big(\Delta \varphi\big) = \lambda_{\xi} \cdot \mathcal{F}\varphi$, for $\varphi \in C_c^{\infty}(G/K)^K$, so, by continuity $\mathcal{F} \big((1- \Delta) f \big) \; = \; (1-\lambda_{\xi}) \cdot \mathcal{F}f$ for all$f \in H^{2\ell +2}(X)$.
\end{proof}
\end{prop}

Let $\mu$ be the multiplication map  $\mu(v) (\xi)\, = \, (1-\lambda_{\xi}) \cdot v(\xi) \, = \, (1 + |\rho|^2+ |\xi|^2) \cdot v(\xi)$ where $\rho$ is the half sum of positive roots.  For $\ell \in \Z$, the weighted $L^2$-spaces $V^{2\ell} \, = \, \{ v \text{ measurable } : \mu^{\ell}(v) \in L^2(\Xi) \}$ with norms
$$\lVert v \rVert_{V^{2\ell}} \;\; = \;\; \lVert \mu^{\ell}(v) \rVert_{L^2(\Xi)} \;\; = \;\; \int_{\Xi} (1 + |\rho|^2+ |\xi|^2)^{\ell} \; |v(\xi)|^2 \; |\mathbf{c}(\xi)|^{-2} \; d\xi$$
are Hilbert spaces with $V^{2\ell + 2} \subset V^{2\ell}$ for all $\ell$.  In fact, these are dense inclusions, since truncations are dense in all $V^{2\ell}$-spaces.  The multiplication map $\mu$ is a Hilbert space isomorphism $\mu: V^{2\ell + 2} \to V^{2\ell}$, since for $v \in V^{2\ell +2}$,
$$\lVert \mu(v) \rVert_{V^{2\ell}} \;\; = \;\; \lVert \mu^{\ell+1}(v) \rVert_{L^2(\Xi)} \;\; = \;\; \lVert v \rVert_{V^{2\ell+2}}$$
The negatively indexed spaces are the Hilbert space duals of their positively indexed counterparts, by integration.  The adjoints to inclusion maps are genuine inclusions, since $V^{2\ell+2} \hookrightarrow V^{2\ell}$ is dense for all $\ell \geq 0$, and, under the identification $(V^{2\ell})^{\ast} = V^{-2\ell}$ the adjoint map $\mu^{\ast}:(V^{2\ell})^{\ast} \to (V^{2\ell +2})^{\ast}$ is the multiplication map $\mu:V^{-2\ell} \to V^{-2\ell -2}$.

\begin{prop}\label{sph_trans_isom}  For $\ell \geq 0$, the spherical transform is an isometric isomorphism $H^{2\ell}(X) \to V^{2\ell}$.

\begin{proof}
On compactly supported functions, the spherical transform $\mathcal{F}$ and its inverse $\mathcal{F}^{-1}$ are given by integrals, which are certainly continuous linear maps.  The Plancherel theorem extends $\mathcal{F}$ and $\mathcal{F}^{-1}$ to isometries between $L^2(X)$ and $L^2(\Xi)$.  Thus $\mathcal{F}$ on $H^{2\ell}(X) \subset L^2(X)$ is a continuous linear $L^2$-isometry onto its image.  

Let $f \in H^{2\ell}(X)$. By Proposition \ref{large_small}, the distributional derivatives $(1-\Delta)^k \, f$ lie in $L^2(X)$ for all $k \leq \ell$.  By the Plancherel theorem and Proposition \ref{diff_components},
$$\lVert (1-\Delta)^{\ell} \, f \rVert_{L^2(X)} \;\; = \;\; \lVert \mathcal{F}\big((1-\Delta)^{\ell} f\big)\rVert_{L^2(\Xi)} \;\;= \;\; \lVert(1- \lambda_{\xi})^{\ell} \cdot \mathcal{F}f \rVert_{L^2(\Xi)}$$
Thus $\mathcal{F}\big(H^{2\ell}(X)\big) \subset V^{2\ell}$.  The following claim shows that $\mathcal{F}^{-1}\big(V^{2\ell}\big) \subset H^{2\ell}(X)$ and finishes the proof. \noqed \end{proof}

\begin{claim*} For $v \in V^{2\ell}$, the distributional derivatives $(1-\Delta)^k \, \mathcal{F}^{-1}v$ lie in $L^2(X)$, for all $0 \leq k \leq \ell$.
\begin{proof}
For test function $\varphi$, the Plancherel theorem implies
$$\big((1-\Delta) \, \mathcal{F}^{-1}v \big) \varphi \;\; = \;\; \mathcal{F}^{-1}v \, \big( (1-\Delta) \varphi \big) \;\; = \;\; v \big( \mathcal{F}\, (1-\Delta) \varphi \big) $$
By Proposition \ref{diff_components} and the Plancherel theorem,
$$ v \big( \mathcal{F}\, (1-\Delta) \varphi \big) \;\; = \;\;  v \big((1-\lambda_{\xi})\cdot  \mathcal{F} \varphi \big) \;\; = \;\; \big((1-\lambda_{\xi})\cdot v\big) \, (\mathcal{F} \varphi) \;\; = \;\; \mathcal{F}^{-1} \big( (1-\lambda_{\xi}) \cdot v \big) \, \varphi$$
By induction, we have the following identity of distributions:  $(1-\Delta)^k \; \mathcal{F}^{-1} v \;\; = \;\;  \mathcal{F}^{-1} \big( (1-\lambda_{\xi})^k \, v \big)$.  Since $\mathcal{F}$ is an $L^2$-isometry and $(1-\lambda_{\xi})^{k} \, v \in L^2(\Xi)$ for all $0 \leq k \leq \ell$, $(1-\Delta)^k \; \mathcal{F}^{-1} v$ lies in $L^2(X)$ for $0 \leq k \leq \ell$.
\end{proof}
\end{claim*}
\end{prop}

\begin{note}\label{spectral_defn}
This Hilbert space isomorphism $\mathcal{F}: H^{2\ell} \to V^{2\ell}$ gives a \emph{spectral} characterization of the $2\ell^{\text{th}}$ Sobolev space, namely the preimage of $V^{2\ell}$ under $\mathcal{F}$.
$$H^{2\ell}(X) = \{ f \in L^2(X) : (1-\lambda_{\xi})^{\ell} \cdot \mathcal{F}f(\xi) \in L^2(\Xi) \}$$
\end{note}

\subsection{Negatively indexed Sobolev spaces and distributions}

Negatively indexed Sobolev spaces allow the use of spectral theory for solving differential equations involving certain \emph{distributions.}

\begin{definition} 
For $\ell >0$, the Sobolev space $H^{-\ell}(X)$ is the Hilbert space dual of $H^{\ell}(X)$.
\end{definition}

\begin{note} Since the space of test functions is a dense subspace of $H^{\ell}(X)$ with $\ell >0$, dualizing gives an inclusion of $H^{-\ell}(X)$ into the space of distributions.  The adjoints of the dense inclusions $H^{\ell} \hookrightarrow H^{\ell -1}$ are inclusions $H^{-\ell+1}(X) \hookrightarrow H^{-\ell}(X)$, and the self-duality of $H^0(X) = L^2(X)$ implies that $H^{\ell}(X) \hookrightarrow H^{\ell - 1}$ for all $\ell \in \Z$.
\end{note}

\begin{prop} \label{sph_trans_isom_neg} The spectral transform extends to an isometric isomorphism on negatively indexed Sobolev spaces $\mathcal{F}:H^{-2\ell} \to V^{-2\ell}$, and for any $u \in H^{2\ell}(X)$, $\ell \in \Z$,  $\mathcal{F} ((1-\Delta) \, u) \, = \, (1-\lambda_{\xi}) \cdot \mathcal{F} u$.

\begin{proof}To simplify notation, for this proof let $H^{2\ell} = H^{2\ell}(X)$.   Propositions \ref{diff_components} and \ref{sph_trans_isom} give the result for positively indexed Sobolev spaces, expressed in the following commutative diagram,
\begin{equation*}
\xymatrix@C=3.5pc@R=3.5pc{
\dots \ar[r]^{(1-\Delta)} & H^4 \ar[r]^{(1-\Delta)} \ar[d]_{\mathcal{F}}^{\approx} & H^2 \ar[r]^{(1-\Delta)} \ar[d]_{\mathcal{F}}^{\approx} & H^0 \ar[d]_{\mathcal{F}}^{\approx}\\
\dots \ar[r]^{\mu} & V^4 \ar[r]^{\mu} & V^2 \ar[r]^{\mu} & V^0 }
\end{equation*}
where $\mu(v) (\xi) \;\; = \;\; (1-\lambda_{\xi}) \cdot v(\xi)$, as above.  Dualizing, we immediately have the commutativity of the adjoint diagram.
\begin{equation*}
\xymatrix@C=3.5pc@R=3.5pc{
(H^0)^{\ast} \ar[r]^{(1-\Delta)^{\ast}} & (H^2)^{\ast} \ar[r]^{(1-\Delta)^{\ast}}  & (H^4)^{\ast} \ar[r]^{(1-\Delta)^{\ast}} & \dots \\
(V^0)^{\ast} \ar[r]^{\mu^{\ast}} \ar[u]_{\approx}^{\mathcal{F}^{\ast}}& (V^{2})^{\ast} \ar[r]^{\mu^{\ast}} \ar[u]_{\approx}^{\mathcal{F}^{\ast}}& (V^{4})^{\ast} \ar[r]^{\mu^{\ast}} \ar[u]_{\approx}^{\mathcal{F}^{\ast}} & \dots }
\end{equation*}
The self-duality of $L^2$ and the Plancherel theorem allow the two diagrams to be connected.
\begin{equation*}
\xymatrix@C=3.5pc@R=3.5pc{
\dots \ar[r]^{(1-\Delta)} & H^4 \ar[r]^{(1-\Delta)} \ar[d]_{\mathcal{F}}^{\approx} & H^2 \ar[r]^{(1-\Delta)} \ar[d]_{\mathcal{F}}^{\approx} & H^0 \ar@/_1pc/[d]_{\mathcal{F}}^{\approx} \ar[r]^{(1-\Delta)^{\ast}} & H^{-2} \ar[r]^{(1-\Delta)^{\ast}}  & H^{-4} \ar[r]^{(1-\Delta)^{\ast}} & \dots\\
\dots \ar[r]^{\mu} & V^4 \ar[r]^{\mu} & V^2 \ar[r]^{\mu} & V^0  \ar[r]^{\mu} \ar@/_1pc/[u]^{\approx}_{\mathcal{F}^{-1}} & V^{-2} \ar[r]^{\mu} \ar[u]_{\approx}^{\mathcal{F}^{\ast}}& V^{-4} \ar[r]^{\mu} \ar[u]_{\approx}^{\mathcal{F}^{\ast}} & \dots }
\end{equation*}
Since $V^{2\ell + 2}$ is dense in $V^{2\ell}$ for all $\ell \in \Z$, and $H^{2\ell} \approx V^{2\ell}$ for all $\ell \in\Z$, $H^{2\ell + 2}$ is dense in $H^{2\ell}$ for all $\ell \in \Z$.  Thus test functions are dense in \emph{all} the Sobolev spaces.  The adjoint map $(1-\Delta)^{\ast}: H^{-2\ell} \to H^{-2\ell -2}$ is the continuous extension of $(1-\Delta)$ from the space of test functions, since, for a test function $\varphi$, identified with an element of $H^{-2\ell}$ by integration,
$$\big((1 - \Delta)^{\ast} \Lambda_{\varphi}\big)(f) \;\; = \;\; \Lambda_{\varphi}\big((1-\Delta)f\big) \;\; = \;\; \langle \varphi, (1-\Delta)f \rangle \;\; = \;\; \langle (1-\Delta) \varphi, f \rangle \;\; = \;\; \Lambda_{(1-\Delta) \varphi} (f)$$
for all $f$ in $H^{2\ell+2}$ by integration by parts, where $\Lambda_{\varphi}$ is the distribution associated with $\varphi$ by integration and $\langle \, , \, \rangle$ denotes the usual inner product on $L^2(G/K)$.  The map $\big(\mathcal{F}^{\ast}\big)^{-1}$ on $H^{-2\ell}$ is the continuous extension of $\mathcal{F}$ from the space of test functions, since for a test function $\varphi$, 
$$\big(\mathcal{F}^{\ast} \, \Lambda_{ \mathcal{F} \varphi}\big) (f ) \;\;= \;\; \Lambda_{\mathcal{F} \varphi}\big( \mathcal{F}f \big) \;\; = \;\; \langle \mathcal{F} \varphi , \mathcal{F} f \rangle_{V^{2\ell}} \;\; = \;\; \langle \varphi, f \rangle_{H^{2\ell}} \;\; = \;\; \Lambda_{\varphi}(f)$$
for all $f \in H^{2\ell}$.  Thus, the following diagram commutes.
\begin{equation*}
\xymatrix@C=3.5pc@R=3.5pc{
\dots \ar[r]^{(1-\Delta)} & H^4 \ar[r]^{(1-\Delta)} \ar[d]_{\mathcal{F}}^{\approx} & H^2 \ar[r]^{(1-\Delta)} \ar[d]_{\mathcal{F}}^{\approx} & H^0 \ar[d]_{\mathcal{F}}^{\approx} \ar[r]^{(1-\Delta)} & H^{-2} \ar[r]^{(1-\Delta)} \ar[d]_{\mathcal{F}}^{\approx} & H^{-4} \ar[r]^{(1-\Delta)} \ar[d]_{\mathcal{F}}^{\approx} & \dots\\
\dots \ar[r]^{\mu} & V^4 \ar[r]^{\mu} & V^2 \ar[r]^{\mu} & V^0  \ar[r]^{\mu} & V^{-2} \ar[r]^{\mu} & V^{-4} \ar[r]^{\mu}  & \dots }
\end{equation*}
In other words, the relation $\mathcal{F} \big((1- \Delta) \, u \big) = (1-\lambda_{\xi}) \cdot \mathcal{F}u$ holds for any $u$ in a Sobolev space.
\end{proof}
\end{prop}

Recall that, for a smooth manifold $M$, the positively indexed \emph{local} Sobolev spaces $H^{\ell}_{\text{loc}}(M)$ consist of functions $f$ on $M$ such that for all points $x \in M$, all open neighborhoods $U$ of $x$ small enough that there is a diffeomorphism $\Phi: U \to \R^n$ with $\Omega = \Phi(U)$ having compact closure, and all test functions $\varphi$ with support in $U$, the function $(f \cdot \varphi) \circ \Phi^{-1} : \Omega \longrightarrow \C$ is in the Euclidean Sobolev space $H^{\ell}(\Omega)$. The Sobolev embedding theorem for local Sobolev spaces states that $H^{\ell +k}_{\text{loc}}(M) \, \subset \, C^k(M)$ for $\ell > \mathrm{dim}(M)/2$.


\begin{prop}\label{global_sob_emb}  For $\ell \, > \, \mathrm{dim}(G/K)/2$,  $H^{\ell + k}(X) \; \subset \; H^{\ell +k}(G/K) \; \subset \; C^k(G/K)$.
\begin{proof}
Since positively indexed global Sobolev spaces on $G/K$ lie inside the corresponding local Sobolev spaces, $H^{\ell}_{\text{loc}}(G/K)  \subset C^k(G/K)$ by local Sobolev embedding.
\end{proof}
\end{prop}

This embedding of global Sobolev spaces into $C^k$-spaces is used to prove that the integral defining spectral inversion for test functions can be extended to sufficiently highly indexed Sobolev spaces, i.e. the abstract isometric isomorphism $\mathcal{F}^{-1} \circ \mathcal{F}: H^{\ell}(X) \to H^{\ell}(X)$ is given by an integral that is convergent uniformly pointwise, when $\ell > \mathrm{dim}(G/K)/2$, as follows. 

 \begin{prop}\label{cvgc_sp_expn}  For $f \in H^{s}(X)$ , $s > k + \mathrm{dim}(G/K)/2$,
$$f \;\; = \;\; \int_{\Xi} \mathcal{F}f(\xi) \, \varphi_{\rho + i \xi} \, |\mathbf{c}(\xi)|^{-2} \, d\xi \;\;\;\;\; \text{in } H^s(X) \text{ and } C^k(X)$$
\begin{proof}
Let $\{\Xi_n\}$ be a nested family of compact sets in $\Xi$ whose union is $\Xi$, $\chi_n$ be the characteristic function of $\Xi_n$, and $f_n$ be given by the following $C^{\infty}(X)$-valued Gelfand-Pettis integral (see \ref{gpints_moll})
$$f_{n} \;\; = \;\; \int_{\Xi} \chi_n(\xi) \, \mathcal{F}f(\xi) \, \varphi_{\rho + i \xi} \; | \mathbf{c}(\xi)|^{-2} \, d\xi$$
Since $ \chi_n(\xi) \, \mathcal{F}f(\xi) $ is compactly supported, $f_n \; = \; \mathcal{F}^{-1} (\chi_n \cdot \mathcal{F}f)$.  Thus, by Propositions \ref{sph_trans_isom} and \ref{sph_trans_isom_neg}
$$\lVert f_n - f_m \rVert_{H^s(X)} \;\; = \;\; \lVert  (\chi_n -\chi_m) \cdot \mathcal{F}f \rVert_{V^{s}}$$
Since $\mathcal{F}f$ lies in $V^{s}$, these tails certainly approach zero as $n, m \to \infty$.  Similarly,
$$\lVert f_n - f \rVert_{H^s(X)} \;\; = \;\; \lVert  (\chi_n -1) \cdot \mathcal{F}f \rVert_{V^{s}} \;\; \longrightarrow \;\; 0 \;\;\;\;\;\;\;\; \text{ as } n \to \infty$$
By Proposition \ref{global_sob_emb}, $f_n$ approaches $f$ in $C^k(X)$.\end{proof}
\end{prop}

The embedding of global Sobolev spaces into $C^k$-spaces also implies that \emph{compactly supported} distributions lie in global Sobolev spaces, as follows.

\begin{prop} \label{cs_distns_sob} Any compactly supported distribution on $X$ lies in a global zonal spherical Sobolev space.  Specifically, a compactly supported distribution of order $k$ lies in $H^{-s}(X)$ for all $s \, > \, k + \mathrm{dim}(G/K)/2$.

\begin{proof}
A compactly supported distribution $u$ lies in $\big(C^{\infty}(G/K)\big)^{\ast}$.  Since compactly supported distributions are of finite order, $u$ extends continuously to $C^k(G/K)$ for some $k \geq 0$.  Using Proposition \ref{global_sob_emb} and dualizing, $u$ lies in $H^{-(\ell+k)}(X)$, for $\ell \, > \, \mathrm{dim}(G/K)/2$.
\end{proof}
\end{prop}

\begin{note} In particular, this implies that the Dirac delta distribution at the base point $x_o = 1 \cdot K$ in $G/K$ lies in $H^{-\ell}(X)$ for all $\ell \, > \, \mathrm{dim}(G/K)/2$.\end{note}

\begin{prop}\label{sph_trans_on_distns}  For a compactly supported distribution $u$ of order $k$,  $\mathcal{F}u \, = \, u(\varphi_{\rho + i \xi})$ in $V^{-s}$ where  $s > k + \mathrm{dim}(G/K)/2$.
\begin{proof}
By Proposition \ref{cs_distns_sob}, a compactly supported distribution $u$ of order $k$ lies in $H^{-s}$ for any $s  >k + \mathrm{dim}(G/K)/2$.  Let $f$ be any element in $H^{s}(X)$.  Then,
$$\langle \mathcal{F}f, \mathcal{F}u \rangle_{V^{s} \times V^{-s}} \;\; = \;\; \langle f, u \rangle_{H^{s}\times V^{s}} \;\; = \;\; u(f)$$
Since the spectral expansion of $f$ converges to it in the $H^{s}(X)$ topology by Proposition \ref{cvgc_sp_expn}, 
$$u(f) \;\; = \;\; u \bigg(\lim_n \; \int_{\Xi_n} \mathcal{F}f(\xi) \, \varphi_{\rho + i\xi} \; |\mathbf{c}(\xi)|^{-2} \, d\xi\bigg) \;\; = \;\; \lim_n \; u \bigg(\int_{\Xi_n} \mathcal{F}f(\xi) \, \varphi_{\rho + i\xi} \; |\mathbf{c}(\xi)|^{-2} \, d\xi\bigg)$$
Since the integral is a $C^{\infty}(X)$-valued Gelfand-Pettis integral (see \ref{gpints_moll}) and $u$ is an element of $(C^{\infty}(X))^{\ast}$,
$$u \bigg(\int_{\Xi_n} \mathcal{F}f(\xi) \, \varphi_{\rho + i\xi} \; |\mathbf{c}(\xi)|^{-2} \, d\xi\bigg) \;\; = \;\; \int_{\Xi_n} \mathcal{F}f(\xi) \, u(\varphi_{\rho + i\xi}) \; |\mathbf{c}(\xi)|^{-2} \, d\xi$$
The limit as $n\to \infty$ is finite, by comparison with the original expression which surely is finite, and thus
$$\langle \mathcal{F}f, \mathcal{F}u \rangle_{V^{s} \times V^{-s}} \;\; = \;\; \int_{\Xi} \mathcal{F}f(\xi) \, u(\varphi_{\rho + i\xi}) \; |\mathbf{c}(\xi)|^{-2} \, d\xi \;\; = \;\; \langle \mathcal{F}f, u(\varphi_{\rho + i \xi}) \rangle_{V^s \times V^{-s}}$$
Thus, $\mathcal{F}u = u(\varphi_{\rho + i \xi})$ as elements of $V^{-s}$.
\end{proof}
\end{prop}

\begin{note} This implies that the spherical transform of the Dirac delta distribution is $\mathcal{F}\delta = \varphi_{\rho + i \xi}(1) = 1$.\end{note}

\subsection{Gelfand-Pettis integrals and mollification} \label{gpints_moll}

We describe the vector-valued (weak) integrals of Gelfand \cite{gelfand36} and Pettis \cite{pettis38} and summarize the key results; see \cite{garrett2011}.  For $X,\mu$ a measure space and $V$ a locally convex, quasi-complete topological vector space, a Gelfand-Pettis (or weak) integral is a vector-valued integral $C_c^o(X,V) \to V$ denoted $f \to I_f$ such that, for all $\alpha \in V^{\ast}$,  $\alpha(I_f) \; = \; \int_X \alpha \circ f \; d\mu$, where this latter integral is the usual scalar-valued Lebesgue integral.  

\begin{note}
 Hilbert, Banach, Frechet, LF spaces, and their weak duals are locally convex, quasi-complete topological vector spaces; see \cite{garrett2011}.
\end{note}

\begin{theorem}\label{gp_ints}
\emph{(i)} Gelfand-Pettis integrals \emph{exist}, are \emph{unique}, and satisfy the following \emph{estimate}:
$$I_f \; \in \; \mu(\mathrm{spt}f) \cdot \big(\text{closure of compact hull of } f(X) \big)$$
\emph{(ii)} Any continuous linear operator between locally convex, quasi-complete topological vetor spaces $T : V \to W$ commutes with the Gelfand-Pettis integral: $T(I_f) \, = \, I_{Tf}$.
\end{theorem}

For a locally compact Hausdorff topological group $G$, with Haar measure $dg$, acting continuously on a locally convex, quasi-complete vector space $V$, the group algebra $C_c^o(G)$ acts on $V$ by averaging:
$$\eta \cdot v \; = \; \int_G \eta(g) \, g \cdot v \, dg$$

\begin{theorem}\label{smooth_vectors_dense} \emph{(i)}
Let $G$ be a locally compact Hausdorff topological group acting continuously on a locally convex, quasi-complete vector space $V$.  Let $\{\psi_i\}$ be an approximate identity on $G$.  Then, for any $v \in V$, $\psi_i \cdot v \to v$ in the topology of $V$.

\emph{(ii)} If $G$ is a Lie group and $\{\eta_i\}$ is a smooth approximate identity on $G$, the mollifications $\eta_i \cdot v$ are smooth.  In particular, for $X \in \mathfrak{g}$,  $X \cdot (\eta \cdot v) \,= \, (L_X \eta) \cdot v$.  Thus the space $V^{\infty}$ of smooth vectors is dense in $V$.
\end{theorem}

\begin{note} For a \emph{function space} $V$, the space of smooth vectors $V^{\infty}$  is \emph{not necessarily} the subspace of \emph{smooth functions} in $V$.  Thus Theorem \ref{smooth_vectors_dense} does \emph{not} prove the density of smooth \emph{functions} in $V$.\end{note}

\section{Fundamental Solution for $(\Delta - \lambda_z)^{\nu}$ on $G/K$}
\label{free_sp_fund_soln}

We now return to prove the main theorem, stated at the beginning of the paper, which gives an explicit expression for the fundamental solution for $(\Delta - \lambda_z)^{\nu}$ on a symmetric space $G/K$.

\subsection{Proof of main theorem}
\label{free_sp_fund_soln_pvthm}

As above, let $G$ be a complex semi-simple Lie group with finite center and $K$ a maximal compact subgroup.  Let $G = NAK$, $\mathfrak{g} = \mathfrak{n} + \mathfrak{a} + \mathfrak{k}$ be corresponding Iwasawa decompositions.  Let $\Sigma$ denote the set of roots of $\mathfrak{g}$ with respect to $\mathfrak{a}$, let $\Sigma^+$ denote the subset of positive roots (for the ordering corresponding to $\mathfrak{n}$), and let $\rho = \tfrac{1}{2} \sum_{\alpha \in \Sigma^+} m_{\alpha} \alpha$, $m_{\alpha}$ denoting the multiplicity of $\alpha$.  Since $G$ is complex, $m_{\alpha} = 2$, for all $\alpha \in \Sigma^+$, so $\rho = \sum_{\alpha \in \Sigma^+} \alpha$.  Let $\mathfrak{a}_{\C}^{\ast}$ denote the set of complex-valued linear functions on $\mathfrak{a}$.  Consider the differential equation on the symmetric space $X = G/K$:
$$(\Delta - \lambda_z)^{\nu} \; u_z \;\; = \;\; \delta_{1 \cdot K}$$
where the Laplacian $\Delta$ is the image of the Casimir operator for $\mathfrak{g}$, $\lambda_z$ is $z^2 - |\rho|^2$ for a complex parameter $z$, $\nu$ is an integer, and $\delta_{1\cdot K}$ is Dirac delta at the basepoint $x_o = 1 \cdot K \in G/K$.  Since $\delta_{1 \cdot K}$ is also \emph{left}-$K$-invariant, we construct a left-$K$-invariant solution on $G/K$ using the harmonic analysis of spherical functions.

\begin{prop}\label{sp_expn_fund_soln}For integral $\nu \, > \, \mathrm{dim}(G/K)/2$, $u_z$ is a continuous left-$K$-invariant function on $G/K$ with the following spectral expansion:
$$u_z(g) \;\; = \;\; \int_{\Xi} \;\frac{(-1)^{\nu}}{(|\xi|^2 + z^2)^{\nu}} \;\; \varphi_{\rho + i \xi}(g) \, |\mathbf{c}(\xi)|^{-2} \, d\xi$$


\begin{proof}
Since $\delta_{1\cdot K}$ is a compactly supported distribution of order zero, by Proposition \ref{cs_distns_sob}, it lies in the global zonal spherical Sobolev spaces $H^{-\ell}(X)$ for all $\ell > \mathrm{dim}(G/K)/2$.  Thus there is an element $u_z$ of $H^{-\ell + 2 \nu}(X)$ satisfying this equation.  The solution $u_z$ is unique in Sobolev spaces, since any $u'_z$ satisfying $\big(\Delta - (z^2 - |\rho|^2)\big)^{\nu} \, u'_z \, = \, \delta_{1 \cdot K}$ must necessarily have the same spherical transform.  For $\nu > \mathrm{dim}(G/K)/2$, the solution is continuous by Proposition \ref{global_sob_emb}, and by Proposition \ref{cvgc_sp_expn}, 
$$u_z(g) \; = \; \int_{\Xi} \mathcal{F}u_z(\xi) \, \varphi_{\rho + i \xi}(g) \, |\mathbf{c}(\xi)|^{-2} \, d\xi\;\; = \;\; \int_{\Xi} \;\frac{(-1)^{\nu}}{(|\xi|^2 + z^2)^{\nu}} \;\; \varphi_{\rho + i \xi}(g) \, |\mathbf{c}(\xi)|^{-2} \, d\xi$$
\end{proof}
\end{prop}

For a \emph{complex} semi-simple Lie group, the zonal spherical functions are \emph{elementary}.  The spherical function associated with the principal series $I_{\chi}$ with $\chi = e^{\rho + i \lambda}$, $\lambda \in \mathfrak{a}_{\C}^{\ast}$ is
$$\varphi_{\rho + i \lambda} \;\; = \;\; \frac{\pi^+(\rho)}{\pi^+(i \lambda)} \; \frac{\sum \mathrm{sgn}(w) \; e^{i \, w \lambda}}{\sum \mathrm{sgn}(w) \, e^{w \rho}}$$
where the sums are taken over the elements $w$ of the Weyl group, and the function $\pi^+$ is the product  $\pi^+(\mu) \; = \; \prod_{\alpha >0} \langle \alpha, \mu \rangle$ over positive roots, \emph{without} multiplicities.  The ratio of $\pi^+(\rho)$ to $\pi^+(i \lambda)$ is the $\mathbf{c}$-function, $\mathbf{c}(\lambda)$.  The denominator can be rewritten
$$\sum_{w \in W} \mathrm{sgn}(w) \, e^{w \rho} \;\;= \;\; \prod_{\alpha \in \Sigma^+} 2 \sinh(\alpha)$$


\begin{prop} \label{int_repn}  The fundamental solution $u_z$ has the following integral representation:
$$u_z \; = \; \frac{(-1)^{\nu} \, (-i)^d}{\pi^+(\rho) \;  \prod 2 \sinh \alpha} \; \cdot \; \int_{\mathfrak{a}^{\ast}} \frac{1}{(|\lambda|^2 + z^2)^{\nu}} \; \pi^+(\lambda) \; e^{i\lambda} \; d\lambda$$
\begin{proof}
In the case of complex semi-simple Lie groups, the inverse spherical transform has an elementary form
\begin{eqnarray*}
\mathcal{F}^{-1} \, f  & = & \int_{\mathfrak{a}^{\ast}/W} f(\lambda) \; \varphi_{\rho + i \lambda} \, |\mathbf{c}(\lambda)|^{-2} \, d\lambda\\
 & = & \int_{\mathfrak{a}^{\ast}/W} f(\lambda) \;  \frac{\pi^+(\rho)}{\pi^+(i \lambda)} \; \frac{\sum \mathrm{sgn}(w) \; e^{i \, w \lambda}}{\sum \mathrm{sgn}(w) \, e^{w \rho}} \; \bigg| \frac{\pi^+(i \lambda)}{\pi^+(\rho)} \bigg|^{2}  \; d\lambda\\
 & = & \frac{1}{\pi^+(\rho) \; \prod 2 \sinh \alpha} \; \int_{\mathfrak{a}^{\ast}/W} f(\lambda) \; \bigg(\sum_{w \in W} \mathrm{sgn}(w) \; e^{i \, w\lambda}\bigg) \; \overline{\pi^+}(i \lambda) \;  d\lambda\end{eqnarray*}
The function $\pi^{+}$ is a homogeneous polynomial of degree $d$, equal to the number of positive roots, counted without multiplicity, so $\overline{\pi^+}(i\lambda) \, = \, (-i)^d \cdot \pi^+(\lambda)$.  Also, $\pi^+$ is $W$-equivariant by the sign character, so the change of variables $\lambda \to w^{-1} \lambda$ yields
 $$\mathcal{F}^{-1} \, f  \;\; = \;\; \frac{(-i)^d }{\pi^+(\rho) \,\prod 2 \sinh \alpha} \cdot \int_{\mathfrak{a}^{\ast}} f(\lambda)\; \pi^+(\lambda) \, e^{i\lambda} \, d\lambda$$
By Proposition \ref{sp_expn_fund_soln}, 
$$u_z \; = \; \frac{(-1)^{\nu} \, (-i)^d}{\pi^+(\rho) \; \prod 2 \sinh \alpha} \; \cdot \; \int_{\mathfrak{a}^{\ast}} \frac{1}{(|\lambda|^2 + z^2)^{\nu}} \; \pi^+(\lambda) \; e^{i\lambda} \; d\lambda$$
\end{proof}
\end{prop}

Let $I(\log a)$ denote the integral we need to compute:
$$I(\log a) \; = \; \int_{\mathfrak{a}^{\ast}} \frac{1}{(|\lambda|^2 + z^2)^{\nu}} \; \pi^+(\lambda) \;  e^{i \langle \lambda, \log a\rangle}  \, d\lambda$$

\begin{prop}\label{red_to_R} The integral $I(\log a)$ can be reduced to an integral over the real line:
$$ I(\log a) \; = \; i^d \, \pi^+(\log a) \cdot \frac{\Gamma(\nu-d)}{\Gamma(\nu)}  \cdot \frac{\pi^{(n-1)/2}}{\Gamma(\nu-d)} \cdot \Gamma\big(\nu-d-\tfrac{n-1}{2}\big) \cdot \int_{\R}  \frac{ e^{i \lambda_1 |\log a| }}{(\lambda_1^2 + z^2)^{\nu-d-(n-1)/2}} \; d\lambda_1 $$
where $n = \mathrm{dim} \, \mathfrak{a}$ and $d$ is the number of positive roots, counted without multiplicity.

\begin{proof} 
Apply the identity
$$\frac{\Gamma(s)}{z^s} \; = \; \int_0^{\infty} t^s \, e^{-tz} \, \frac{dt}{t}$$
to $(|\lambda|^2 + z^2)^{-\nu}$ in the integrand of $I(\log a)$:
\begin{eqnarray*}
I(\log a) & = & \frac{1}{\Gamma(\nu)} \cdot \int_0^{\infty} \; \int_{\mathfrak{a}^{\ast}} t^{\nu} \; e^{-t(|\lambda|^2 + z^2)} \; \pi^+(\lambda) \;  e^{i \lambda}  \, d\lambda \, \frac{dt}{t}\\
 & = & \frac{1}{\Gamma(\nu)} \cdot \int_0^{\infty} t^{\nu} \;e^{-t z^2}\; \int_{\mathfrak{a}^{\ast}}  e^{-t |\lambda|^2} \; \pi^+(\lambda) \;  e^{i \langle \lambda, \log a\rangle}  \, d\lambda  \, \frac{dt}{t}
\end{eqnarray*}
Change variables $\lambda \to \lambda/\sqrt{t}$.
\begin{eqnarray*}
I(\log a) & = & \frac{1}{\Gamma({\nu})} \cdot \int_0^{\infty} t^{\nu} \;e^{-t z^2}\; \int_{\mathfrak{a}^{\ast}}  e^{-|\lambda|^2} \; t^{-d/2} \cdot \pi^+(\lambda) \;  e^{i \langle \lambda/\sqrt{t} \,, \, \log a\rangle}  \;\; t^{-n/2} \,d\lambda \, \frac{dt}{t}\\
 & = & \frac{1}{\Gamma(\nu)} \cdot \int_0^{\infty} t^{\nu - (d + n)/2} \;e^{-t z^2}\; \int_{\mathfrak{a}^{\ast}}  e^{-|\lambda|^2} \;  \pi^+(\lambda) \;  e^{-i \langle \lambda \, , \, -\log a/\sqrt{t}\rangle}  \;\; d\lambda \, \frac{dt}{t}\\
\end{eqnarray*}
The polynomial $\pi^+$ is in fact \emph{harmonic}.  See, for example, Lemma 2 in \cite{urakawa} or, for a more direct proof, Theorem \ref{pi_plus_harmonic}, below.  Thus the integral over $\mathfrak{a^{\ast}}$ is the Fourier transform of the product of a Gaussian and a harmonic polynomial, and by Hecke's identity,
\begin{eqnarray*}
\int_{\mathfrak{a}^{\ast}}  e^{-|\lambda|^2} \;  \pi^+(\lambda) \;  e^{-i \langle \lambda \, , \, -\log a/\sqrt{t}\rangle}  \;\; d\lambda  & = & (-i)^d \, \pi^+ \big(-\log a/\sqrt{t} \big) \, e^{-|\log a|^2/t} \\
 & = & i^d \, t^{-d/2} \, \pi^+(\log a) \, e^{-|\log a|^2/t}
\end{eqnarray*}
Returning to the main integral, 
\begin{eqnarray*}
I(\log a)  & = & \frac{ i^d \, \pi^+(\log a)}{\Gamma(\nu)}  \cdot \int_0^{\infty} t^{\nu - d - n/2} \; e^{-t z^2} \, e^{-|\log a|^2/t} \, \frac{dt}{t}\\
 & = & \frac{ i^d \, \pi^+(\log a)}{\Gamma(\nu)}  \cdot \int_0^{\infty} t^{\nu - d} \; e^{-t z^2} \, \big(t^{-n/2} e^{-|\log a|^2/t} \big)\, \frac{dt}{t}
\end{eqnarray*}
Replacing the Gaussian by its Fourier transform, 
\begin{eqnarray*}
I(\log a) & = & \frac{ i^d \, \pi^+(\log a)}{\Gamma(\nu)}  \cdot \int_0^{\infty} t^{\nu - d} \; e^{-t z^2} \,  \int_{\mathfrak{a}^{\ast}} e^{i \langle \lambda, \log a \rangle} \; e^{-t|\lambda|^2} \; d\lambda \, \frac{dt}{t}\\
 & = & \frac{ i^d \, \pi^+(\log a)}{\Gamma(\nu)}  \cdot \int_{\mathfrak{a}^{\ast}} \; \int_0^{\infty}  t^{\nu - d} \; e^{-t (|\lambda|^2 +z^2)}  \, \frac{dt}{t} \;\; e^{i \langle \lambda, \log a \rangle} \; d\lambda\\
  & = & \frac{ i^d \, \pi^+(\log a)}{\Gamma(\nu)}  \cdot \int_{\mathfrak{a}^{\ast}} \; \frac{\Gamma(\nu-d)}{(|\lambda|^2 +z^2)^{\nu-d}} \; e^{i \langle \lambda, \log a \rangle} \; d\lambda\\
 & = & i^d \, \pi^+(\log a) \cdot \frac{\Gamma(\nu-d)}{\Gamma(\nu)}  \cdot \int_{\mathfrak{a}^{\ast}} \; \frac{e^{i \langle \lambda, \log a \rangle}}{(|\lambda|^2 +z^2)^{\nu-d}} \; d\lambda\\
\end{eqnarray*}
We denote this integral by $J(|\log a|)$, since it is is rotation-invariant as a function of $\log a$.  Writing $\lambda = (\lambda_1, \dots , \lambda_n)$, we may assume $\langle \lambda, \log a\rangle = \lambda_1 \cdot |\log a|$, and then
$$J(|\log a|) \; = \;  \int_{\mathfrak{a}^{\ast}} \; \frac{e^{i \lambda_1 \cdot |\log a|}}{(|\lambda|^2 +z^2)^{\nu-d}} \; d\lambda$$
Identifying $\mathfrak{a}^{\ast}$ with $\R^n$, $J(|\log a|)$ is
\begin{eqnarray*}
\lefteqn{\frac{1}{\Gamma(\nu-d)} \cdot \int_0^{\infty} t^{\nu-d} \; \int_{\R^n} e^{-t(|\lambda|^2 + z^2)} \;e^{i \lambda_1 |\log a| } \; d\lambda_1 \, d\lambda_2 \, \dots \, \lambda_n \, \frac{dt}{t}}\\
 & = & \frac{1}{\Gamma(\nu-d)} \cdot \int_0^{\infty} t^{\nu-d} \int_{\R} \ e^{-t(\lambda_1^2+ z^2)} \, e^{i \lambda_1 |\log a| }\; d\lambda_1 \;\; \int_{\R^{n-1}} e^{-t(\lambda_2^2 +\dots +\lambda_n^2)} \, d\lambda_2 \, \dots \, \lambda_n \; \frac{dt}{t} \\
 & = & \frac{\pi^{(n-1)/2}}{\Gamma(\nu-d)} \cdot \int_0^{\infty} t^{\nu-d -(n-1)/2} \,  \int_{\R} e^{-t(\lambda_1^2+z^2)} \; e^{i \lambda_1 |\log a| } \; d\lambda_1 \, \frac{dt}{t} \\
& = &  \frac{\pi^{(n-1)/2}}{\Gamma(\nu-d)} \cdot \Gamma\big(\nu-d-\tfrac{n-1}{2}\big) \cdot \int_{\R}  \frac{ e^{i \lambda_1 |\log a| }}{(\lambda_1^2 + z^2)^{\nu-d-(n-1)/2}} \; d\lambda_1
\end{eqnarray*}
Thus we have the desired conclusion, since
$$I(\log a) \; = \;  i^d \, \pi^+(\log a) \cdot \frac{\Gamma(\nu-d)}{\Gamma(\nu)}  \cdot J(|\log a|) $$
\end{proof}
\end{prop}

The remaining integral can be evaluated by \emph{residues} when the exponent in the denominator is a sufficiently large \emph{integer}, i.e. when $G$ is of $\emph{odd}$ rank.  For the even rank case, the integral can be expressed in terms of a \emph{K-Bessel function}.

\begin{prop}\label{odd_rk} For $n$ odd, let $\nu \; = \; d + \frac{n+1}{2}$.   Then
$$I(\log a) \; = \; i^d \, \pi^+(\log a) \cdot \frac{\pi^{(n+1)/2}}{\Gamma(d + (n+1)/2)} \; \cdot \; \frac{e^{-z \, |\log a|}}{z}$$
\begin{proof}
Recall from Proposition \ref{red_to_R} that $I(\log a)$ is
$$ i^d \, \pi^+(\log a) \cdot \frac{\Gamma(\nu-d)}{\Gamma(\nu)}  \cdot \frac{\pi^{(n-1)/2}}{\Gamma(\nu-d)} \cdot \Gamma\big(\nu-d-\tfrac{n-1}{2}\big) \cdot \int_{\R}  \frac{ e^{i \lambda_1 |\log a| }}{(\lambda_1^2 + z^2)^{\nu-d-(n-1)/2}} \; d\lambda_1$$
Taking $\nu \; = \; d + \frac{n+1}{2}$,
$$\int_{\R}  \frac{ e^{i \lambda_1 |\log a| }}{(\lambda_1^2 + z^2)^{\nu-d-(n-1)/2}} \; d\lambda_1 \; = \; \int_{\R}  \frac{ e^{i \lambda_1 |\log a| }}{\lambda_1^2 + z^2} \; d\lambda_1 \;\; = \;\; 2\pi i \cdot \operatorname*{Res}_{\lambda_1 = iz} \; \bigg( \frac{ e^{i \lambda_1 |\log a| }}{\lambda_1^2 + z^2} \bigg) \;\; = \;\;  \frac{\pi \, e^{-z \, |\log a|}}{z}$$
Thus,
$$I(\log a)  \;\; = \;\; i^d \, \pi^+(\log a) \cdot \frac{\pi^{(n+1)/2}}{\Gamma(d + (n+1)/2)} \; \cdot \; \frac{e^{-z \, |\log a|}}{z}$$
\end{proof}
\end{prop}

\begin{prop}\label{even_rk} For $n$ even, let $\nu \; = \; d + \frac{n}{2} + 1$.  Then,
$$I(\log a) \;\;= \;\; i^d \, \pi^+(\log a) \cdot  \frac{\Gamma((n/2)+1)}{\Gamma(d + (n/2) + 1)} \; \cdot \;\;  \frac{\pi^{n/2}}{\Gamma((n/2) + 1) } \cdot \frac{ |\log a|}{z} \cdot  K_1(z \, |\log a|)$$
\begin{proof}
Recall from Proposition \ref{red_to_R} that
$$I(\log a) \;\; = \;\;i^d \, \pi^+(\log a) \cdot \frac{\Gamma(\nu-d)}{\Gamma(\nu)}  \cdot \frac{\pi^{(n-1)/2}}{\Gamma(\nu-d)} \cdot \Gamma\big(\nu-d-\tfrac{n-1}{2}\big) \cdot \int_{\R}  \frac{ e^{i \lambda_1 |\log a| }}{(\lambda_1^2 + z^2)^{\nu-d-(n-1)/2}} \; d\lambda_1 $$
The integral over $\R$ is
\begin{eqnarray*}
\int_{\R}  \frac{ e^{i \lambda_1 |\log a| }}{(\lambda_1^2 + z^2)^{3/2}} \; d\lambda_1  & = & \int_{\R} \frac{\cos(\lambda_1 |\log a|)}{(\lambda_1^2 + z^2)^{3/2}} \, d\lambda_1 + \int_{\R} \frac{i\sin(\lambda_1 |\log a|)}{(\lambda_1^2 + z^2)^{3/2}} \, d\lambda_1\\
 & = & 2 \cdot \int_{0}^{\infty} \frac{\cos(\lambda_1 |\log a|)}{(\lambda_1^2 + z^2)^{3/2}} \, d\lambda_1\\
 &=& \frac{\sqrt{\pi} \, |\log a|}{\Gamma(3/2) \, z} \; \cdot \;  K_{1}(z\, |\log a|)
\end{eqnarray*}
where $K_1$ is a modified Bessel function $K_{\alpha}$ of the second kind, which has the following integral representation (see \cite{abramowitz-stegun1964}, 9.6.25).
$$K_{\alpha}(xz) \;\; = \;\; \frac{\Gamma(\alpha+ \frac{1}{2}) \, (2z)^{\alpha}}{\pi^{1/2} \, x^{\alpha}} \cdot \int_0^{\infty} \frac{\cos(xt)}{(t^2 + z^2)^{\alpha + 1/2}} \, dt \;\;\;\;\; \mathrm{Re}(\alpha) > -\tfrac{1}{2}, \; x >0, \; |\mathrm{arg} \, z| < \tfrac{\pi}{2}$$
Thus,
$$I(|\log a|)\;\;= \;\; i^d \, \pi^+(\log a) \cdot \frac{\Gamma((n/2)+1)}{\Gamma(d + (n/2) + 1)} \; \cdot \;  \frac{\pi^{n/2}}{\Gamma((n/2) + 1) } \cdot \frac{ |\log a|}{z} \cdot  K_1(z \, |\log a|)$$
\end{proof}
\end{prop}

Now we prove the theorem stated in the introduction: 

\begin{theorem}\label{fmla_for_fundl_soln}  When $G$ is of odd rank, let $\nu = d + \frac{n+1}{2}$, where $d$ is the number of positive roots, counted without multiplicities, and $n$ is the rank.  Then the fundamental solution $u_z$ for the operator $(\Delta - \lambda_z)^{\nu}$ on $G/K$ is given by:
$$u_z(a) \;\; = \;\;  \frac{(-1)^{d+(n+1)/2} \, \pi^{(n+1)/2}}{\pi^+(\rho)\, \Gamma(d+(n+1)/2) } \; \cdot \; \prod_{\alpha\in \Sigma^+}  \, \frac{\alpha(\log a)}{2 \sinh(\alpha(\log a))} \; \cdot \;\; \frac{e^{-z|\log a|}}{z}  $$
When $G$ is of even rank, let $\nu  =  d + \frac{n}{2} + 1$.  Then, with $K_1$ the usual Bessel function,
$$u_z(a) \;\; = \;\; \frac{(-1)^{d + (n/2) + 1} \,\pi^{n/2}}{\pi^+(\rho) \; \Gamma(d + (n/2) + 1)} \; \cdot  \; \prod_{\alpha\in \Sigma^+}  \, \frac{\alpha(\log a)}{2 \sinh(\alpha(\log a))} \;\cdot \; \frac{ |\log a|}{ z} \; \cdot \; K_1(z \, |\log a|)$$
\begin{proof}
Recall that the fundamental solution is
$$u_z(a) \;\; = \;\; \frac{(-1)^{\nu} \, (-i)^d}{\pi^+(\rho) \; \prod 2 \sinh(\alpha(\log a))} \; \cdot \; I(\log a)$$
By Proposition \ref{odd_rk}, when $G$ is of odd rank and $\nu = d + \frac{n+1}{2}$, 
\begin{eqnarray*}
u_z(a) & = &  \frac{(-1)^{d+(n+1)/2} \, (-i)^d }{\pi^+(\rho) \; \prod 2 \sinh(\alpha(\log a))} \; \cdot  \; i^d \, \pi^+(\log a) \cdot \frac{\pi^{(n+1)/2}}{\Gamma(d + (n+1)/2)} \; \cdot \; \frac{e^{-z \, |\log a|}}{z}
\end{eqnarray*}
\begin{eqnarray*}
 & = &   \frac{(-1)^{d+(n+1)/2} \, \pi^{(n+1)/2}}{\pi^+(\rho)\, \Gamma(d+(n+1)/2) } \; \cdot \; \,\frac{ \pi^+(\log a)}{\prod 2 \sinh(\alpha(\log a))}  \cdot \;\; \frac{e^{-z|\log a|}}{z}  
\end{eqnarray*}
\begin{eqnarray*}
 & = &   \frac{(-1)^{d+(n+1)/2} \, \pi^{(n+1)/2}}{\pi^+(\rho)\, \Gamma(d+(n+1)/2) } \; \cdot \; \,\prod_{\alpha \in \Sigma^+} \frac{\alpha(\log a)}{2 \sinh(\alpha(\log a))}  \cdot \;\; \frac{e^{-z|\log a|}}{z}  
\end{eqnarray*}
By Proposition \ref{even_rk}, when $G$ is of even rank and $\nu= d + \frac{n}{2} + 1$,
\begin{eqnarray*}
 u_z(a) & = & \frac{(-1)^{d + (n/2) + 1}}{\pi^+(\rho)} \; \cdot \; \frac{ \pi^+(\log a)}{\prod 2 \sinh(\alpha(\log a))} \cdot \frac{\Gamma((n/2)+1)}{\Gamma(d + (n/2) + 1)} \; \cdot \; J(|\log a|) 
 \end{eqnarray*}
 \begin{eqnarray*}
 & = &  \frac{(-1)^{d + (n/2) + 1} \,\pi^{n/2}}{\pi^+(\rho) \; \Gamma(d + (n/2) + 1)} \; \cdot \; \prod_{\alpha \in \Sigma^+} \frac{\alpha(\log a)}{2 \sinh(\alpha(\log a))} \; \cdot \; \frac{ |\log a|}{ z} \; \cdot \; K_1(z \, |\log a|)
\end{eqnarray*}
\end{proof}
\end{theorem}

\begin{note}\label{asymptotic} For fixed $\alpha$, large $|z|$, and $\mu = 4 \alpha^2$ (see \cite{abramowitz-stegun1964}, 9.7.2),
$$K_{\alpha}(z) \approx \sqrt{\tfrac{\pi}{2z}}   e^{-z} \; \big( 1   +  \frac{\mu-1}{8z} + \frac{(\mu-1)(\mu-9)}{2! \, (8z)^2} + \frac{(\mu-1)(\mu -9)(\mu -25)}{3! \, (8z)^3}  + \dots \big)\;\; \big(|\mathrm{arg} \, z | < \tfrac{3\pi}{2}\big)$$ 
Thus in the even rank case the fundamental solution has the following asymptotic:
$$u_z(a) \; \approx \; \frac{(-1)^{d + (n/2) + 1} \,\pi^{(n+1)/2}}{\sqrt{2} \, \pi^+(\rho) \; \Gamma(d + (n/2) + 1)} \; \cdot \; \prod_{\alpha \in \Sigma^+} \frac{\alpha(\log a)}{2 \sinh(\alpha(\log a))}\; \cdot \; \sqrt{\frac{ |\log a|}{ z}} \; \cdot \; \frac{e^{-z |\log a|}}{z}$$
\end{note}


%

\begin{note} For any integer $\nu > \mathrm{dim}(G/K)/2 = n + d/2$, the argument used in the proof of Proposition \ref{even_rk} gives a formula for the bi-$K$-invariant fundamental solution for $(\Delta - \lambda_z)^{\nu}$ in terms of K-Bessel functions:
$$u_z(a) \;\; = \;\; \frac{2 (-1)^{\nu}}{\pi^+(\rho)\Gamma(\nu)} \; \cdot \prod_{\alpha \in \Sigma^+} \frac{\alpha(\log a)}{2 \sinh(\log a)} \cdot \; \bigg(\frac{|\log a|}{2z}\bigg)^{\nu - d - n/2} \; \cdot \; K_{\nu -d -n/2}(z \, |\log a|)$$
In the odd rank case, the classical formula for Bessel functions of half integer order
$$K_{m + 1/2}(z) \;\; = \;\; \sqrt{\frac{\pi}{2}} \; \cdot \; \frac{e^{-z}}{\sqrt{z}} \; \cdot \; \sum_{j = 0}^m \; \frac{(m+j)!}{j!(m-j)!} \; (2z)^{-j}$$
allows us to write the fundamental solution in elementary terms.  Let $\nu = m + d + (n+1)/2$, where $m$ is any non-negative integer.  Then
$$u_z(a) \;\; = \;\; \frac{(-1)^{m + d + (n+1)/2} \, \pi^{(n+1)/2}}{(m + d + (n-1)/2)! \, \pi^+(\rho)} \; \cdot \; \prod_{\alpha \in \Sigma^+} \frac{\alpha(\log a)}{2 \sinh(\alpha(\log a))} \; \cdot \; \frac{e^{-z |\log a|}}{z} \; \cdot \; P(|\log a|, z^{-1})$$
where $P$ is a degree $m$ polynomial in $|\log a|$ and a degree $2m$ polynomial in $z^{-1}$.

\end{note}

\begin{note} Recall from Proposition \ref{sp_expn_fund_soln} that zonal spherical Sobolev theory ensures the continuity of $u_z$ for $\nu$ chosen as in the theorem.  For $G = SL_2(\C)$, the continuity is visible, since fundamental solution is, up to a constant,
$$u_z(a_r) \;\; = \;\; \frac{r \, e^{-(2z-1)r}}{(2z-1) \, \sinh r} \;\;\;\;\; \text{where } \; a_r \;= \; \left(\begin{smallmatrix} e^{r/2} & 0 \\ 0 & e^{-r/2} \end{smallmatrix}\right)$$
\end{note}

\begin{note} As Brian Hall has pointed out, there is a simple relation between Laplacian on $G/K$ and the Laplacian on $\mathfrak{p}$ when $G$ is complex, which allows us to determine the fundamental solution on $G/K$ in terms of the Euclidean fundamental solution on $\mathfrak{p}$ (\cite{hall-mitchell2005}, proof of Theorem 2).  (See also Helgason's discussion of the wave equation on $G/K$ in \cite{helgason1994}.)  Let $J$ be the function on $\mathfrak{p}$ defined by
$$\int_{G/K} f(g) \, dg \;\; = \;\; \int_{\mathfrak{p}} f(\exp(X)) \, J(X) \, dX$$
Then, for $\lambda_z = z^2 - |\rho|^2$, as above,
$$(\Delta_{G/K} - \lambda_z)^{\nu} \, f \circ \exp \;\; = \;\; J^{-1/2} \; \big(\Delta_{\mathfrak{p}} - z^2\big)^{\nu} \, J^{1/2} f \circ \exp$$
When $G$ is complex, the restriction to $\mathfrak{a}$ has an elementary square root, 
$$J^{1/2}(H) \; = \;\; \prod_{\alpha \in \Sigma^+} \frac{\sinh(\alpha(H))}{\alpha(H)}$$
Indeed the integral we evaluated after applying Hecke's identity is an integral representation for the Euclidean fundamental solution on $\mathfrak{a}^{\ast} \approx \R^n$.
\end{note}

\subsection{The harmonic property of $\pi^{+}$}
\label{free_sp_fund_soln_piplusharm}

Let G be complex semi-simple Lie group.  We will give a direct proof that the function $\pi^+: \mathfrak{a}^{\ast} \to \R$ given by  $\pi^+(\mu) \; = \; \prod_{\alpha >0} \langle \alpha, \mu \rangle$ where the product is taken over all postive roots, counted without multiplicity, is harmonic with respect to the Laplacian naturally associated to the pairing on $\mathfrak{a}^{\ast}$.  (See also \cite{urakawa}, Lemma 2, where this result is obtained as a simple corollary of the less trivial fact that $\pi^+$ divides any polynomial that is $W$-equivariant by the sign character.)  It is this property that enables us to use Hecke's identity in the computations above.  We will use the following lemma.  
\begin{lemma}\label{sum_zero}
Let $I$ be the set of all non-orthogonal pairs of distinct positive roots, as functions on $\mathfrak{a}$.  Then $\pi^+$ is harmonic if $ \sum_{(\beta, \gamma) \, \in \, I} \frac{\langle \beta, \gamma \rangle}{ \beta\,  \gamma } \; = \; 0$.
\begin{proof}
Considering $\mathfrak{a}^{\ast}$ as a Euclidean space, its Lie algebra can be identified with itself.  For any basis $\{ x_i \}$ of $\mathfrak{a}^{\ast}$, the Casimir operator (Laplacian) is $\Delta \; = \; \sum_{i} x_i \, x_i^{\ast}$.  For any $\alpha, \beta$ in $\mathfrak{a}^{\ast}$ and any $\lambda \in \mathfrak{a}$
\begin{eqnarray*}
\lefteqn{\Delta \; \langle \alpha, \lambda \rangle  \langle \beta, \lambda \rangle \;\; = \;\;   \sum_{i} x_i \,  \big(\langle \alpha, x_i^{\ast} \rangle  \, \langle \beta , \lambda \rangle +   \langle \alpha, \lambda \rangle \,  \langle \beta , x_i^{\ast} \rangle \big)}\\
& = &  \sum_{i}  \big(\langle \alpha, x_i^{\ast} \rangle  \, \langle \beta , x_i  \rangle +   \langle \alpha, x_i  \rangle \,  \langle \beta , x_i^{\ast} \rangle \big) \;\; = \;\; 2 \langle \alpha, \beta \rangle
\end{eqnarray*}
Thus
\begin{eqnarray*}
\lefteqn{\Delta \pi^+ \;\; = \;\; \sum_i x_i x_i^{\ast} \pi^+ \;\; = \;\; \sum_i x_i \sum_{\beta >0} \alpha(x_i^{\ast}) \cdot \frac{\pi^+}{\beta}}\\
& = &  \sum_i \sum_{\beta >0} \beta(x_i^{\ast})  \cdot \bigg( \sum_{\gamma \neq \beta} \gamma(x_i) \cdot \frac{\pi^+}{\beta \gamma} \bigg) \; = \; \bigg( \sum_{\beta \neq \gamma} \frac{\langle \beta, \gamma \rangle}{\beta \, \gamma } \bigg) \; \cdot \pi^+
\end{eqnarray*}
and $\pi^+$ is harmonic if the sum in the statement of the Lemma is zero.
\end{proof}
\end{lemma}

\begin{note} \label{semisimple} When the Lie algebra $\mathfrak{g}$ is not simple, but merely semi-simple, i.e. $\mathfrak{g} = \mathfrak{g}_1 \oplus \mathfrak{g}_2$, any pair $\beta, \gamma$ of roots with $\beta \in \mathfrak{g}_1$ and $\gamma \in \mathfrak{g}_2$ will have $\langle \beta,\gamma\rangle = 0$, so it suffices to consider $\mathfrak{g}$ \emph{simple}.\end{note}

\begin{prop} \label{rk_two}
For $\mathfrak{g} = \mathfrak{sl_3}$, $\mathfrak{sp}_2$, or $\mathfrak{g}_2$, the following sum over all pairs $(\beta, \gamma)$ of distinct positive roots is zero: $ \sum_{\beta \neq \gamma} \frac{\langle \beta, \gamma \rangle}{ \beta\,  \gamma } \; = \; 0$.
\begin{proof}
The positive roots in $\mathfrak{sl}_3$ are $\alpha$, $\beta$, and $(\alpha + \beta)$ with  $\langle \alpha, \alpha \rangle \; = \; 2$, $\langle \beta, \beta \rangle \; = \; 2$, $\langle \alpha, \beta \rangle = -1$.  In other words, the two simple roots have the same length and have an angle of $2\pi/3$ between them.  The pairs of distinct positive roots are $(\alpha, \beta)$, $(\alpha, \alpha + \beta)$ and $(\beta, \alpha + \beta)$, so the sum to compute is
$$\frac{\langle \alpha, \beta \rangle }{\alpha \, \beta } \;  + \; \frac{\langle \alpha, \alpha + \beta \rangle }{ \alpha \, (\alpha + \beta) } \; + \; \frac{\langle \beta, \alpha + \beta \rangle}{\beta \, (\alpha + \beta)}$$
Clearing denominators and evaluating the parings,
$$\langle \alpha, \beta \rangle \cdot (\alpha + \beta) \;  + \; \langle \alpha, \alpha + \beta \rangle\cdot \beta  \; + \; \langle \beta, \alpha + \beta \rangle \cdot \alpha \;\;\; = \;\;\; - (\alpha + \beta) \, + \, \beta \, + \, \alpha \;\;\; = \;\;\; 0$$
For $\mathfrak{sp}_2$, the simple roots have lengths 1 and $\sqrt{2}$ and have an angle of $3\pi/4$ between them: $\langle \alpha, \alpha \rangle \; = \; 1$, $\langle \beta, \beta \rangle \; = \; 2$, $\langle \alpha, \beta \rangle \; = \; -1$.  The other positive roots are $(\alpha + \beta)$ and $(2\alpha + \beta)$.  The non-orthogonal pairs of distinct positive roots are $(\alpha, \beta)$, $(\alpha, 2 \alpha + \beta)$, $(\beta, \alpha + \beta)$, and $(\alpha + \beta, 2\alpha + \beta)$.  So the sum we must compute is
$$\frac{\langle \alpha, \beta \rangle}{\alpha \, \beta} \; + \; \frac{\langle \alpha, 2\alpha + \beta \rangle}{\alpha \, (2\alpha + \beta)} \; + \; \frac{\langle \beta, \alpha + \beta \rangle}{\beta \, (\alpha + \beta)} \; + \; \frac{\langle \alpha + \beta , 2\alpha + \beta \rangle}{(\alpha + \beta)(2\alpha + \beta)} $$
Again, clearing denominators,
$$\langle \alpha, \beta \rangle \cdot (\alpha + \beta) (2\alpha + \beta) \; + \; \langle \alpha, 2\alpha + \beta \rangle  \cdot\beta (\alpha + \beta) \; + \; \langle \beta, \alpha + \beta \rangle \cdot \alpha (2\alpha + \beta) \; + \; \langle \alpha + \beta , 2\alpha + \beta \rangle  \cdot  \alpha \beta $$
and evaluating the pairings,
\begin{eqnarray*}
\lefteqn{ -(\alpha + \beta) (2\alpha + \beta) \; + \; \beta (\alpha + \beta) \; + \; \alpha (2\alpha + \beta) \; + \;  \alpha \beta}\\
& = &   -(2\alpha^2 + 3\alpha \beta + \beta^2) \; + \; \alpha\beta + \beta^2 \; + \; 2\alpha^2 +\alpha\beta \; + \;  \alpha \beta \;\;\; = \;\;\; 0
\end{eqnarray*}
Finally we consider the exceptional Lie algebra $\mathfrak{g}_2$.  The simple roots have lengths 1 and $\sqrt{3}$ and have an angle of $5\pi/6$ between them:  $\langle \alpha, \alpha \rangle = 1$, $\langle \beta, \beta \rangle \; = \; 3$, $\langle \alpha, \beta \rangle \; = \; -3/2$.  The other positive roots are $ (\alpha + \beta)$, $(2\alpha + \beta)$, $(3\alpha + \beta)$, and $(3\alpha + 2\beta)$.  Notice that the roots $\alpha$ and $\alpha + \beta$ have the same length and have an angle of $3\pi/2$ between them.  So together with their sum $2\alpha + \beta$, they form a copy of the $\mathfrak{sl}_3$ root system.  The three terms corresponding to the three pairs of roots among these roots will cancel, as in the $\mathfrak{sl}_3$ case.  Similarly, the roots $(3\alpha + \beta)$ and $\beta$ have the same length and have an angle of $3\pi/2$ between them, so, together with their sum, $(3\alpha + 2\beta)$ they form a copy of the $\mathfrak{sl}_3$ root system, and the three terms in the sum corresponding to the three pairs among these roots will also cancel.  The remaining six pairs of distinct, non-orthogonal postitive roots are $(\alpha, 3\alpha + \beta)$, $(\alpha, \beta)$, $(3\alpha + \beta, 2\alpha + \beta)$, $(2\alpha + \beta, 3\alpha + 2\beta)$, $(3\alpha + 2\beta, \alpha + \beta)$, and $(\alpha + \beta, \beta)$.  We shall see that the six terms corresponding to these pairs cancel as a group.  After clearing denominators, the relevant sum is
\begin{eqnarray*}
\lefteqn{\langle \alpha, \beta \rangle \cdot  (\alpha+\beta)(2\alpha+\beta)(3\alpha + \beta)(3\alpha+2\beta) \;+ \; \langle \alpha, 3\alpha + \beta \rangle \cdot \beta (\alpha+\beta)(2\alpha + \beta)(3\alpha+2\beta)}\\
&+ & \langle 3\alpha + \beta, 2\alpha + \beta  \rangle \cdot \alpha \beta(\alpha+\beta)(3\alpha+2\beta) +  \langle 2\alpha + \beta, 3\alpha + 2\beta \rangle \cdot \alpha \beta (\alpha+\beta)(3\alpha + \beta)\\
& + & \langle 3\alpha + 2\beta, \alpha + \beta \rangle \cdot \alpha \beta(2\alpha + \beta)(3\alpha+\beta) +  \langle \alpha + \beta, \beta \rangle \cdot \alpha (2\alpha + \beta)(3\alpha +\beta)(3\alpha + 2\beta)
\end{eqnarray*}
Evaluating the pairings and factoring out $(3/2)$, this is
\begin{eqnarray*}
\lefteqn{-(\alpha+\beta)(2\alpha+\beta)(3\alpha + \beta)(3\alpha+2\beta)\; + \; \beta (\alpha+\beta)(2\alpha + \beta)(3\alpha+2\beta)}\\
& + & \alpha \beta(\alpha+\beta)(3\alpha+2\beta) \; + \; \alpha \beta (\alpha+\beta)(3\alpha + \beta)\\
& + & \alpha \beta(2\alpha + \beta)(3\alpha+\beta) \;+ \;  \alpha (2\alpha + \beta)(3\alpha +\beta)(3\alpha + 2\beta)
\end{eqnarray*}
Multiplying out,
\[
\begin{array}{rcrcrcrcr}
 -18 \alpha^4 & - & 45 \, \alpha^3 \beta &- &  40\, \alpha^2 \beta ^2 & - & 15\, \alpha \beta^3 & - & 2 \beta^4\\
& + & 6 \, \alpha^3\beta & + & 13\, \alpha^2\beta^2 &+& 9\, \alpha \beta^3 &+& 2\beta^4 \\
& + & 3\,\alpha^3 \beta &+& 5 \,\alpha^2 \beta^2  &+& 2\, \alpha \beta^3\\
& + & 3 \, \alpha ^3 \beta  &+& 4 \,\alpha^2 \beta^2  &+& \alpha \beta^3\\
& + & 6 \,\alpha^3 \beta  &+& 5\, \alpha^2 \beta^2  &+& \alpha \beta^3 \\
+ 18\alpha^4 & + & 27 \, \alpha^3 \beta  &+& 13 \, \alpha^2 \beta^2  &+& 2 \,\alpha \beta^3\\
\end{array}
\]
This sum is zero.
\end{proof}
\end{prop}

\begin{prop}\label{arb_rk}
For any complex simple Lie algebra $\mathfrak{g}$, the following sum over all pairs $(\beta, \gamma)$ of distinct positive roots is zero: $ \sum_{\beta \neq \gamma} \frac{\langle \beta, \gamma \rangle}{ \beta\,  \gamma } \; = \; 0$.

\begin{proof}
Let $I$ be the indexing set $\{ (\beta, \gamma) \}$ of pairs of distinct, non-orthogonal positive roots.  For each $(\beta, \gamma) \in I$, let $\mathcal{R}_{\beta, \gamma}$ be the two-dimensional root system generated by $\beta$ and $\gamma$.  For such a root system $\mathcal{R}$, let $I_{\mathcal{R}}$ be the set of pairs of distict, non-orthogonal positive roots, where positivity is inherited from the ambient $\mathfrak{g}$.  The collection $J$ of all such $I_{\mathcal{R}}$ is a cover of $I$.  We refine $J$ to a subcover $J'$ of disjoint sets, in the following way.  

For any pair $I_{\mathcal{R}}$ and $ I_{\mathcal{R}'}$ of sets in $J$ with non-empty intersection, there is a two-dimensional root system $\mathcal{R}''$ such that $I_{\mathcal{R''}}$ contains $I_{\mathcal{R}}$ and $I_{\mathcal{R}'}$.  Indeed, letting $(\beta,\gamma)$ and $(\beta', \gamma')$ be pairs in $I$ generating $\mathcal{R}$ and $\mathcal{R}'$ respectively, the non-empty intersection of $I_{\mathcal{R}}$ and $I_{\mathcal{R}'}$ implies that there is a pair $(\beta'', \gamma'')$ lying in both $I_{\mathcal{R}}$ and $I_{\mathcal{R}'}$.  Since $\mathcal{R}$ and $\mathcal{R}'$ are two-dimensional and $\beta''$ and $\gamma''$ are linearly independent, all six roots lie in a plane.  Since all six roots lie in the root system for $\mathfrak{g}$, they generate a two-dimensional root system $\mathcal{R}''$ containing $\mathcal{R}$ and $\mathcal{R}'$, and $I_{\mathcal{R}''} \supset I_{\mathcal{R}}, I_{\mathcal{R}'}$.  Thus we refine $J$ to a subcover $J'$: if $I_{\mathcal{R}}$ in $J$ intersects any $I_{\mathcal{R'}}$ in $J$, replace $I_{\mathcal{R}}$ and $I_{\mathcal{R'}}$ with the set $I_{\mathcal{R''}}$ described above.    The sets $I_{\mathcal{R}}$ in $J'$ are mutually disjoint, and, for any $(\beta, \gamma) \in I$, there is a root system $\mathcal{R}$ such that $(\beta, \gamma) \in I_{\mathcal{R}} \in J'$, thus
$$ \sum_{(\beta, \gamma) \, \in \, I} \frac{\langle \beta, \gamma \rangle}{ \beta\,  \gamma } \; = \; \sum_{I_{\mathcal{R}} \, \in \, J'} \; \sum_{(\beta, \gamma) \, \in \, I_{\mathcal{R}}} \frac{\langle \beta, \gamma \rangle}{ \beta\,  \gamma }$$
By the classification of complex simple Lie algebras of rank two, $\mathcal{R}$ is isomorphic to the root system of $\mathfrak{sl}_3$, $\mathfrak{sp}_2$, or $\mathfrak{g}_2$.  Thus, by Proposition \ref{rk_two}, the inner sum over $I_{\mathcal{R}}$ is zero, proving that the whole sum is zero.  

Note that the refinement \emph{is} necessary, as there are copies of $\mathfrak{sl}_3$ inside $\mathfrak{g}_2$.  Note also that the only time the root system of $\mathfrak{g}_2$ appears is in the case of $\mathfrak{g}_2$ itself, since, by the classification, $\mathfrak{g}_2$ is the only root system containing roots that have an angle of $\pi/6$ or $5\pi/6$ between them.
\end{proof}
\end{prop}

\begin{note} See \cite{hall-stenzel2003}, Lemma 2, for a proof of Proposition \ref{arb_rk} when $G$ is not necessarily complex. \end{note}

\begin{theorem}\label{pi_plus_harmonic} For a complex semi-simple Lie group $G$, the function $\pi^+: \mathfrak{a}^{\ast} \to \R$ given by  $\pi^+(\mu) \; = \; \prod_{\alpha >0} \langle \alpha, \mu \rangle$ where the product is taken over all postive roots, counted without multiplicity, is harmonic with respect to the Laplacian naturally associated to the pairing on $\mathfrak{a}^{\ast}$. 
\begin{proof} This follows immediately from Lemma \ref{sum_zero}, Remark \ref{semisimple}, and Proposition \ref{arb_rk}. \end{proof}
\end{theorem}

\bibliography{fund_solns.bib,harm_an_sph_fcns.bib,gelfand_pettis.bib,sobolev_refs}{}
\bibliographystyle{plain}

\end{document}